\documentclass[12pt,reqno]{amsart}
    \usepackage{amssymb,amsmath,amsthm,newlfont}
    \usepackage{psfrag}
    \usepackage[dvips,final]{graphicx}
\usepackage[T1]{fontenc}
\usepackage{color}

\numberwithin{equation}{section}

\newcommand{\N}{\mathbb{N}}

\newcommand{\C}{\mathbb{C}}

\newcommand{\D}{\mathbb{D}}

\newcommand{\M}{\mathbb{M}}
\newcommand{\Z}{\mathbb{Z}}

\newcommand{\cF}{{\mathcal{F}}}

\newcommand{\cD}{{\mathcal{D}}}

\newcommand{\re}{\mathrm{Re}}

\newcommand{\beqa}{\begin{eqnarray*}}
\newcommand{\bea}{\begin{eqnarray}}
\newcommand{\eeqa}{\end{eqnarray*}}
\newcommand{\eea}{\end{eqnarray}}

\newtheorem{thm}{\bf Theorem}[section]

\newtheorem*{prop*}{\bf Proposition}

\newtheorem{lem}[thm]{\bf Lemma}

\begin{document}

\title{On zero sets in Fock spaces}

\date{\today}

\thanks{ Research partially supported by "Hassan II Academy of Science and Technology".}

\author[D. Aadi, B. Bouya, Y. Omari]{Driss Aadi, Brahim Bouya, Youssef Omari}

\address{ Laboratory of Mathematical Analysis and Applications (LAMA), CeReMAR, Faculty of Sciences,
Mohammed V University in Rabat, 4 Av. Ibn Battouta, Morocco.
}

\email{aadidriss@gmail.com,brahimbouya@gmail.com,omariysf@gmail.com, brahimbouya@fsr.ac.ma}

\begin{abstract}
{We prove that zero sets for distinct Fock spaces are not the same, this is an answer of a question asked by K. Zhu in \cite[Page. 209]{Zhu}.}
\end{abstract}
\maketitle
\section{Introduction and statement of main results}

For $\alpha>0$ and $p>0$ the Fock space $\cF^p_\alpha$ consists of those entire functions $f$ satisfying
$$\|f\|^{p}_{p,\alpha}:=\int_{\C}|f(z)|^p dA_{p\alpha/2}(z)\ <\infty,$$
where $$dA_\beta(z):=\frac{\beta}{\pi}e^{-\beta|z|^2}dA(z),\qquad \beta>0,$$
and $A$ represents the Lebesgue area measure on the complex plane $\C.$
It is known that the space $\cF_\alpha^p$ endowed with the norm $\|\cdot \|_{p,\alpha}$ is a Banach space when $ p\geq1$, while for $p<1$ it is a complete metric space,  see for instance \cite[Chap. 2]{Zhu}.

A sequence $\Lambda$ of complex numbers is called a zero set for $\cF_\alpha^p$ if there exists a function $f\in\cF_\alpha^p\setminus\{0\}$ such that the zero set $\{z\in\C\ :\ f(z)=0\}$ of $f,$ counting multiplicities, coincides with $\Lambda$.  At the present time there is no complete characterization of zero sets for Fock spaces. In  \cite{ZhuZer} and \cite[Chap. 5]{Zhu} K. Zhu has presented many properties enjoyed by zero sets in $\cF_\alpha^p,$ in particular he proved that the spaces $\cF_\alpha^p$  and $\cF_\beta^q$ always possess  different zero sets in the case where  $\alpha\neq\beta,$ regardless of $p$ and $q.$ He then asked whether this remains true if $\alpha=\beta,$ see \cite[Page. 209]{Zhu}.
In this paper, we answer positively to this question by considering a special translation of simple lattices with a uniform density.

The upper and lower Beurling-Landau density of a sequence $\Lambda\subset\C$ is known respectively as the following
$$\cD^+(\Lambda):=\limsup_{\rho\mapsto\infty}\sup_{z\in\C}\frac{N_{\Lambda}(z,\rho)}{\pi\rho^2}$$
and
$$\cD^-(\Lambda):=\liminf_{\rho\mapsto\infty}\inf_{z\in\C}\frac{N_{\Lambda}(z,\rho)}{\pi\rho^2},$$
where $N_{\Lambda}(z,\rho)$ is the number of the elements in the intersection of $\Lambda$ and the Euclidian open disk $\D(z,\rho)$ of center $z\in\C$ and radius $\rho>0$. Our main result is the following theorem.

\begin{thm}\label{thmA}
Let $p$ and $q$ be two positive numbers such that
$p>q.$ There exists a sequence  $\Lambda$ in $\C$ satisfying
\begin{equation}\label{densCond}
\cD^+(\Lambda)=\cD^-(\Lambda)=\alpha/\pi,
\end{equation}
and such that $\Lambda$ is a zero set for $\cF_\alpha^p$ but it is not for $\cF_\alpha^q.$
\end{thm}

The condition \eqref{densCond} in Theorem \ref{thmA} shows that our result is not based on the characterization of sampling and interpolating sets for Fock spaces, given by  K. Seip and R. Wallst\'en \cite{Sei, SW}. Indeed, a sequence of the critical density $\alpha/\pi$ is neither sampling nor interpolating for $\cF^p_\alpha.$ 

\section{Proof of Theorem \ref{thmA}}

We start this section with some well known preliminaries. We consider the following square lattice
$$\Lambda:=\{z_{m,n}:=a(m+in)\ :\ m,n\in\Z\},$$
where $a$ is a positive number and $\Z$ denotes the usual set of integers. The imaginary axis is clearly a line of symmetry for  $\Lambda.$ 
By translating the positive real points of $\Lambda$ away from $0$ and keeping this symmetry unchanged, we define the following modified lattice
$$\Lambda_R:=\{w_{m,n}\ :\ m,n\in\Z\},$$ where $R$ is a positive number and
\begin{equation}\label{gsh}
w_{m,n}:=\left\{
  \begin{array}{lll}
   z_{m,n}, & \hbox{if }\  n\neq0 \text{ or } m=n=0, \\
    a(m+Rm/|m|), & \hbox{if }\  n=0\text{ and }m\neq0.
  \end{array}
\right.
\end{equation}
We observe that if  $R$ is a positive integer, then  $\Lambda_R$ is actually obtained from $\Lambda$ by just removing the following finite symmetric set $\{\pm am\ :\ m\in\{1,2,...,R\}\}.$ The well known Weierstrass function associated to $\Lambda$ is defined
by
$$\sigma_{a}(z):=z\prod_{(m,n)\in\Z^2 \atop (m,n)\neq(0,0)}
\left(1-\frac{z}{z_{m,n}}\right)\exp\left(\frac{z}{z_{m,n}}+\frac{z^2}{2z_{m,n}^2}\right),\qquad z\in\C,$$
one can see the textbooks \cite{Boa,Lev}. The modified Weierstrass  function associated to $\Lambda_R,$ introduced by K. Seip, is given by
$$\sigma_{a,R}(z):=z\prod_{(m,n)\in\Z^2 \atop (m,n)\neq(0,0)}
\left(1-\frac{z}{w_{m,n}}\right)\exp\left(\frac{z}{w_{m,n}}+\frac{z^2}{2z_{m,n}^2}\right),\qquad z\in\C,$$
see for instance \cite[Chap. 4]{Zhu}. For $h_1$ and $h_2$ being two positive functions, we use the following notation  $h_1\lesssim h_2$ to mean that
$h_1\leq c h_2$ for some positive constant $c$. We also write
$h_1\asymp h_2$ if both  $h_1\lesssim h_2$ and $h_2\lesssim h_1.$  In the next section we give the proof of the following lemma.

\begin{lem}\label{lemP} Let $\alpha$ be a positive number. We have
\begin{equation}\label{jouj2}
\left|\sigma_{a,R}(z)\right|e^{-\frac{\alpha}{2}|z|^2}\asymp \frac{d(z,\Lambda_R)}{(1+|z|)^{2R}}, \qquad z\in\C,
\end{equation}
where $a=\sqrt{\pi/\alpha}$ and 
$R$ is a positive constant. 
\end{lem}

We now let $p$ and $q$ be two positive numbers such that
$p>q.$ We take a number $R$ satisfying
 $\displaystyle\frac{1}{p}<R<\frac{1}{q}.$
By using \eqref{jouj2},
\begin{eqnarray}
\left|\sigma_{a,R}(z)\right|^pe^{-\frac{p\alpha}{2}|z|^2}|z|\lesssim
1/|z|^{2pR-1}, \qquad |z|\geq1.
\end{eqnarray}
Since $2pR-1>1,$ we then obtain $\sigma_{a,R}\in\cF_\alpha^p,$ and hence $\Lambda_R$ is a zero set for $\cF_\alpha^p.$ A standard argument  by contradiction shows that $\Lambda_{R}$ cannot be a zero set for $\cF_\alpha^q.$ For the sake of completeness, we sketch here the proof. We suppose that there exists a function
$f\in\cF_\alpha^q\setminus\{0\}$ with zero set $\Lambda_{R}.$ By Hadamard's factorization theorem, we have
$$f(z)=z e^{Q(z)}\prod_{(m,n)\in\Z^2 \atop (m,n)\neq(0,0)}
\left(1-\frac{z}{w_{m,n}}\right)\exp\left(\frac{z}{w_{m,n}}+\frac{z^2}{2w_{m,n}^2}\right), \qquad z\in\C,$$
where $Q$ is a polynomial of degree at most $2.$
We also obviously have
$$\sum_{m\geq1}\frac{R(2m+R)}{m^2(m+R)^2}=M_R<+\infty.$$
Thus, for $z\in\C,$
\begin{eqnarray}\nonumber
f(z)&=& e^{Q(z)}\sigma_{a,R}(z)\prod_{m\neq0}\exp\left(\big(\frac{1}{2w_{m}^2}-\frac{1}{2z_{m}^2}\big)z^2\right)
\\\nonumber&=& \sigma_{a,R}(z)\exp\left(Q(z)-\frac{M_R}{a^2}z^2\right)
\\\label{Hfactor}&=& \sigma_{a,R}(z) e^{L(z)},
\end{eqnarray}
where $L$ is a polynomial of degree at most $2.$
Using \eqref{jouj2} and \eqref{Hfactor}
\begin{eqnarray}\nonumber
\|f\|_{q}^{q}&=& \frac{\alpha q}{2\pi}\int_{\C}|f(z)|^q ~~e^{-\frac{\alpha q}{2}|z|^2} dA(z)
\\\nonumber&=&\frac{\alpha q}{2\pi}\int_{\C}|\sigma^{q}_{a,R}(z)|~~e^{-\frac{\alpha q}{2}|z|^2}|~~e^{qL(z)}|  dA(z)
\\\nonumber&\gtrsim&\int_{|z|\geq\frac{a}{8}}d^q(z,\Lambda_R)|z|^{-2Rq}|~~e^{qL(z)}|  dA(z)
\\\label{estbelow1}&\gtrsim&\int_{\C\setminus\D(\Lambda_R,a/8)}|z|^{-2Rq}~~|e^{qL(z)}|dA(z),
\end{eqnarray}
where $$\D(\Lambda_R,a/8):=\bigcup_{\lambda\in\Lambda_R}\D(\lambda,a/8).$$
For $\lambda\in\Lambda_R\setminus\{0\}$ and a point $w\in\D(\lambda,a/8),$ the subharmonicity of  the function $z\mapsto\phi(z):=|z|^{-2Rq}|~~e^{qL(z)}|$ in $\C\setminus\{0\}$ gives
\begin{eqnarray*}
\phi(w)\lesssim \int_{a/4\leq|z-w|\leq 3a/8}\phi(z) dA(z)
\leq\int_{\D(\lambda,a/2)\setminus\D(\lambda,a/8)}\phi(z) dA(z).
\end{eqnarray*}
It follows
\begin{eqnarray*}
\int_{\D(\lambda,a/8)}\phi(w)dA(w)\lesssim
\int_{\D(\lambda,a/2)\setminus\D(\lambda,a/8)}\phi(z)dA(z).
\end{eqnarray*}
Therefore
\begin{eqnarray}\label{estbelow2}
\int_{\D(\Lambda_R,a/8)\setminus\D(0,a/8)}\phi(w)dA(w)\lesssim
\int_{\C\setminus\D(\Lambda_R,a/8)}\phi(z) dA(z),
\end{eqnarray}
since $$\D(\lambda_1,a/2)\cap\D(\lambda_2,a/2)=\emptyset, \qquad \lambda_1\neq\lambda_2.$$
By combining \eqref{estbelow1} and \eqref{estbelow2}
\begin{eqnarray}\label{estbelow3}
\int_{\C\setminus\D(0,a/8)}\phi(z)dA(z)\lesssim\|f\|_q^q.
\end{eqnarray}
Using again the subharmonicity of $\phi$ and taking account of \eqref{estbelow3}, we deduce that $\phi$ is bounded at $\infty,$ and hence
$z\mapsto e^{qL(z)}$ possesses a polynomial growth. Thus
$L$ is a constant and by consequence
\begin{eqnarray}\label{estbelow4}
\int_{a/8}^{\infty}|z|^{-2Rq+1}d|z|\lesssim\|f\|_q^q.
\end{eqnarray}
The inequality \eqref{estbelow4} is in contradiction with the fact that $2Rq-1<1.$ Hence  $\Lambda_R$ is not a zero set for $\cF_\alpha^q,$ which finishes the proof of Theorem \ref{thmA}.

\section{Proof of Lemma \ref{lemP} \label{prlemA}}
Let $\alpha$ be a positive number and consider the lattice $\Lambda$ generated by $a=\sqrt{\pi/\alpha}.$ By using the symmetry with respect to the imaginary axis enjoyed by the lattices $\Lambda$ and $\Lambda_R,$ we simply compute
\begin{eqnarray}\nonumber
\sigma_{a,R}(z)&=&\sigma_{a}(z)\prod_{m\neq0}\frac{1-z/w_m}{1-z/z_m}\times \frac{\exp(z/w_m)}{\exp(z/z_m)}
\\\label{nisba}&=&\sigma_{a}(z)\prod_{m\geq1}\frac{1-\left(z/w_m\right)^2}{1-\left(z/z_m\right)^2}, \qquad z\in\C\setminus{a\Z},
\end{eqnarray}
where $w_{m}:=w_{m,0}$ and $z_{m}:=z_{m,0}.$ For proving Lemma \ref{lemP}, we claim that it is sufficient to show
\begin{equation}\label{claim}
\psi_R(z):=\prod_{m\geq1}\Big|\frac{m+R-z}{m-z}\Big|\frac{m}{m+R}\asymp \frac{d(z,\Z^+_R)}{d(z,\Z^+)(1+|z|)^{R}},\qquad z\in\C\setminus{\Z^+},
\end{equation}
where  $\Z^+$ is the set of non-negative integers and 
$$\Z^+_R:=\left\{m+R,\quad m\in\Z^+\right\}.$$ Indeed, assume that \eqref{claim} holds. Then
\begin{eqnarray*}
\psi_R(z/a)= \prod_{m\geq 1}\left|\frac{1-\left(z/a(m+R)\right)}{1-\left(z/am\right)}\right| \asymp \frac{d(z,a\Z^+_R)}{d(z,a\Z^+)(1+|z|)^{R}},\qquad z\in\C\setminus{a\Z^+},
\end{eqnarray*}
and
\begin{eqnarray*}
\psi_R(-z/a)= \prod_{m\geq 1}\left|\frac{1+\left(z/a(m+R)\right)}{1+\left(z/am\right)}\right| \asymp\frac{d(z,a\Z^-_R)}{d(z,a\Z^{-})(1+|z|)^{R}},\quad z\in\C\setminus{a\Z^-},
\end{eqnarray*}
where $\Z^-:=-\Z^+$ and $\Z^-_R:=-\Z^+_R.$
We clearly have
$$ \frac{d(z,a\Z^+_R)}{d(z,a\Z^+)}\times\frac{d(z,a\Z^-_R)}{d(z,a\Z^-)}\asymp\frac{d(z,\Lambda_{R})}{d(z,\Lambda)},\qquad z\in\C\setminus{a\Z}. $$
Thus
\begin{eqnarray}\label{c1}
 \prod_{m\geq 1}\left|\frac{1-\left(z/a(m+R)\right)^2}{1-\left(z/am\right)^2}\right| \asymp \frac{d(z,\Lambda_{R})}{d(z,\Lambda)(1+|z|)^{2R}},\qquad z\in\C\setminus{a\Z}.
\end{eqnarray}
By using \eqref{nisba} and \eqref{c1} we deduce
\begin{eqnarray}\label{wahed}
\left|
\sigma_{a,R}(z)\right|\asymp
\frac{|\sigma_{a}(z)|d(z,\Lambda_{R})}{d(z,\Lambda)(1+|z|)^{2R}}, \qquad z\in\C\setminus{a\Z}.
\end{eqnarray}
On the other hand, it is known that 
\begin{equation}\label{jouj}
\left|\sigma_a(z)\right|e^{-\frac{\alpha}{2}|z|^2}\asymp d(z,\Lambda), \qquad z\in\C,
\end{equation}
see for instance \cite[Corollary 1.21]{Zhu}. Hence
\begin{eqnarray}\label{wahed}
\left|
\sigma_{a,R}(z)\right|e^{-\frac{\alpha}{2}|z|^2}\asymp\frac{d(z,\Lambda_{R})}{(1+|z|)^{2R}}, \qquad z\in\C,
\end{eqnarray}
which proves Lemma \ref{lemP}.

Let us now prove \eqref{claim}. For this aim, it is sufficient to consider only the situation when $[R],$ the integer part of $R,$ equals zero. Indeed, we fix  a number $R>1.$ We can factorize $\psi_R$ as follows
\begin{eqnarray*}
\psi_R(z)&=&\psi_\beta(z-[R])\prod_{m=1}^{[R]}\frac{m+\beta}{|m-z|}
\\&\asymp&\psi_\beta(z-[R])\prod_{m=1}^{[R]}\frac{1}{|m-z|},\qquad z\in\C\setminus{\Z^+},
\end{eqnarray*}
where $\beta:=R-[R].$ We have
$$\prod_{m=1}^{[R]}\frac{1}{|m-z|}\asymp \frac{1}{d(z,\{1,2,\cdots,[R]\})(1+|z|)^{[R]-1}},\qquad z\in\C\setminus\{1,2,\cdots,[R]\}.$$
If we show
\begin{equation*}
\psi_\beta(z)\asymp
\frac{d(z,\Z^+_\beta)}{d(z,\Z^+)(1+|z|)^{\beta}},\qquad z\in\C\setminus{\Z^+},
\end{equation*} 
then
\begin{eqnarray*}
\psi_\beta(z-[R])&\asymp&
\frac{d(z,\Z^+_R)}{d(z,\Z^+_{[R]})(1+|z|)^{\beta}},\qquad z\in\C\setminus{\Z^+_{[R]}}.
\end{eqnarray*}
Therefore
\begin{eqnarray*}
\psi_R(z)&\asymp&
\frac{d(z,\Z^+_R)}{d(z,\Z^+)(1+|z|)^{R}},\qquad z\in\C\setminus{\Z^+},
\end{eqnarray*}
which proves \eqref{claim}. So, in the sequel we suppose that $[R]=0.$  We now set
$$\N_{z}:=\{m_z-2,m_z-1,m_z\}\cap\Z^+,$$
where
$$m_{z}:=\min\{m\in\Z^+\ :\ m-x\geq0\},$$
and $x=:\re(z)$ is the real part of $z.$
Since $\N_{z}$ contains at most three elements, we obviously get
$$\prod_{m\in\N_z}\frac{m}{m+R}\asymp1, \qquad z\in\C.$$
If $x\leq1$ we then obtain $m_{z}=1,$ $d(z,\Z^+)=|1-z|$ and $d(z,\Z^{+}_{R})=|1+R-z|,$   and if $x>1$ then $m_{z}\geq2,$ $d(z,\Z^+)=\min\{|m_z-z|,|m_z-1-z|\}=d(z,\N_z)$ and
$d(z,\Z^+_{R})=d(z-R,\N_z).$ Thus
\begin{eqnarray}\label{nihaya1}
\prod_{m\in\N_z}\left|\frac{m+R-z}{m-z}\right|\frac{m}{m+R}\asymp 
\frac{d(z,\Z^+_R)}{d(z,\Z^+)},\qquad z\in\C\setminus{\Z^+}.
\end{eqnarray}
Taking account of \eqref{nihaya1}, for proving \eqref{claim} it remains to show
\begin{eqnarray*}
\prod_{m\in\Z^+\setminus\N_z}\left|\frac{m+R-z}{m-z}\right|\frac{m}{m+R}\asymp (1+|z|)^{-R},\qquad z\in\C,
\end{eqnarray*}
for which it is necessary and sufficient to show that
\begin{eqnarray}\label{nihaya3}
\prod_{m\in\Z^+\setminus\N_z}\left|\frac{m+R-z}{m-z}\right|\frac{m}{m+R}\asymp (1+|z|)^{-R},\qquad  |z|\rightarrow\infty.
\end{eqnarray}
We set
$$\varphi_1(z):=\prod_{ m> 2|z|}\Big|\frac{m+R-z}{m-z}\Big|\frac{m}{m+R},\qquad z\in\C.$$ 
If $m> 2|z|$ then
$|m+R-z|\geq |m-z|> m/2,$ and by using the following usual inequality
$$\log(1+u)\leq u,\qquad u\geq0,$$
we compute
\begin{eqnarray*}
|\log\varphi_1(z)| 
& \leq & \sum_{m>2|z|}\left|\log\Big|\frac{1-z/(m+R)}{1-z/m}\Big| \right|
\\& =& \sum_{m>2|z|}\max\left\{\log\Big|\frac{1-z/(m+R)}{1-z/m}\Big|,\ \log\Big|\frac{1-z/m}{1-z/(m+R)}\Big|\right\}
\\& \leq & \sum_{m>2|z|}\max\left\{ \Big|\frac{Rz}{(m+R)(m-z)}\Big|,\ 
\Big|\frac{Rz}{m(m+R-z)}\Big|
\right\}
\\   & \leq &2R |z|\underset{m> 2|z|}{\sum}\frac{1}{m^{2}}
=O(1).
\end{eqnarray*}
We then deduce
\begin{equation}\label{bona1}
\varphi_1(z)\asymp 1, \qquad z\in\C.
\end{equation}
With \eqref{bona1} in mind, for proving  \eqref{nihaya3} it remains now to show
\begin{equation}\label{lba9ya}
\varphi_2(z):=\prod_{m\in\M_z}\left|\frac{m+R-z}{m-z}\right|\frac{m}{m+R}\asymp (1+|z|)^{-R}, \qquad  |z|\rightarrow\infty,
\end{equation}
where  $$\M_z:=\{m\in\Z^+\setminus\N_z\ :\ m\leq2|z|\}.$$ 
We recall the following classical equality
\begin{equation}\label{logeq}
\log|1+u|=\re(u)+O(|u|^2),\qquad  \re(u)\geq -\frac{1}{2}\text{ and }|u|\leq1.
\end{equation}
By using \eqref{logeq}, 
\begin{eqnarray*}
\sum_{m\in\M_z}\log\frac{m}{m+R} & = & -\sum_{m\in\M_z}\log\left(1+\frac{R}{m}\right) \nonumber\\
    & = & -\sum_{m\in\M_z} \frac{R}{m} + O(1)\nonumber \\
    & = & -R\log(|z|) + O(1), \nonumber
\end{eqnarray*}
which gives
\begin{equation}\label{first}
\prod_{m\in\M_z}\frac{m}{m+R}\asymp|z|^{-R}\asymp (1+|z|)^{-R},\qquad   |z|\rightarrow\infty.
\end{equation}
We have 
$$\sum_{m\in\M_z}
\left|\frac{R}{m-z}\right|^2\leq 
\sum_{m\in\M_z}
\frac{1}{|m-x|^2}
\leq 2\sum_{m\geq1}\frac{1}{m^2}<\infty,$$
and since $$-\frac{1}{2}\leq\re(\frac{R}{m-z})\ \text{ and }\ \frac{R}{|m-z|}\leq 1,\qquad m\in\M_z,$$
then, by using again \eqref{logeq},
\begin{eqnarray}\label{ce1}
\sum_{m\in\M_z}\log \left|1+\frac{R}{m-z}\right|
=R\sum_{m\in\M_z} \frac{m-x}{(m-x)^2+y^2}+O(1),
\end{eqnarray}
where $y$ is the imaginary part of $z.$
For $|z|\geq3/2,$ 
$$\M_{z}^{+}:=\{m\in\M_z\ :\  m-x\geq0\}=\{m_z+1,m_z+2,\cdots,[2|z|]\}\neq\emptyset,$$
and by a simple calculation
\begin{eqnarray}\nonumber\label{cee2}
\sum_{m\in\M_{z}^{+}} \frac{m-x}{(m-x)^2+y^2}
&=&\frac{1}{2}\log\frac{([2|z|]-x)^2+y^2}{(m_z+1-x)^2+y^2}+O(1)
\\&=&\displaystyle\left\{
  \begin{array}{lll}
  O(1), & \hbox{if }\   x\leq0,\\
\log\frac{|z|}{1+|y|}+O(1), & \hbox{if }\  x>0.
  \end{array}
\right.
\end{eqnarray}
We need to distinguish between two different cases.
In the case where $x\leq3,$ we obtain $m_z\leq3$ and hence  $\M_{z}^{+}=\M_z.$ In this case, we either have  $x\leq0$ or $|z|\asymp|y|,$ for $|z|\geq3/2.$  In both situations we deduce the desired estimate \eqref{lba9ya} by combining the estimates \eqref{first}, \eqref{ce1} and \eqref{cee2}.
In the case where $x>3,$ we obtain $m_z\geq4$ and by consequence 
$$\M_z\setminus\M_{z}^{+}=\{1,2,\cdots,m_z-3\}\neq\emptyset.$$ In this case, 
\begin{eqnarray}\nonumber\label{cee1}
\sum_{m\in\M_z\setminus\M_{z}^{+}} \frac{m-x}{(m-x)^2+y^2}
&=&\frac{1}{2}\log\frac{(m_z-3-x)^2+y^2}{(1-x)^2+y^2}+O(1)
\\&=&
\log\frac{1+|y|}{|z|}+O(1).
\end{eqnarray}
We again deduce \eqref{lba9ya} by joining  together \eqref{first}, \eqref{ce1}, \eqref{cee2} and
\eqref{cee1}. The proof of Lemma \ref{lemP} is  completed.

\textbf{Acknowledgement.} The authors would like to thank the reviewer for his/her comments, which helped to improve the manuscript.


\bigskip

\begin{thebibliography}{100}

\bibitem{Boa}
{\sc R. P. Boas}, {\it Entire Functions.} Academic Press, New York, 1954.


\bibitem{Lev}
{\sc B. Y. Levin}, {\it Lectures on Entire Functions.} Transl. Math. Monographs 150, American Mathematical Society, Providence, RI, 1996.

\bibitem{Sei}
{\sc K. Seip}, {\it Density theorems for sampling and interpolation in the Bargmann-Fock space.} Bull. Amer. Math. Soc. 26, 322-328, 1992.

\bibitem{SW}
{\sc K. Seip, R. Wallst\'en}, {\it  
Density theorems for sampling and interpolation in the Bargmann-Fock
space II.} J. Reine Angew. Math. 429, 107-113, 1992.

\bibitem{ZhuZer}
{\sc K. Zhu}, {\it 
Zeros of functions in Fock spaces.} Complex Variables 21, 87-98, 1993.

\bibitem{Zhu}
{\sc K. Zhu}, {\it Analysis on Fock Spaces.}  Springer-Verlag, New York, 2012.


\end{thebibliography}
\end{document}